\documentclass{article}%
\usepackage{amsmath}%
\usepackage{amsfonts}%
\usepackage{amssymb}%
\usepackage{graphicx}

\begin{document}

\title{Functional Power Series}
\author{Henrik Stenlund\thanks{The author is grateful to Visilab Signal Technologies for supporting this work.}
\\Visilab Signal Technologies Oy, Finland}
\date{April 24, 2012}
\maketitle

\begin{abstract}
This work introduces a new functional series for expanding an analytic function in terms of an arbitrary analytic function. It is generally applicable and straightforward to use. It is also suitable for approximating the behavior of a function with a few terms. A new expression is presented for the composite function's n'th derivative. The inverse-composite method is handled in this work also. \footnote{Visilab Report \#2012-04}
\subsection{Keywords}
functional series, function series, function power series of functions, function polynomial series, composite differentiation
\subsection{Mathematical Classification}
Mathematics Subject Classification 2010: 11L03, 30E10, 30K05, 30B10, 30C10, 30D10, 40H05
\end{abstract}
\subsection{}
Dedicated to my father who believed in me.
\section{Introduction}
\subsection{General}
Often the need arises to form a power series of a known function in terms of some other known function, called a functional series. Merely a few methods exist for this purpose, for example the Puiseux series. It expands in terms of fractional powers of the argument. Schl\"{o}milch series is done in terms of zero'th order Bessel function of the first kind. They offer no generality to this problem. The well-known Taylor's power series \cite{Taylor1715} is expanded in terms of a polynomial power series.

The literary background for this subject is thin. Textbooks (\cite{Bender1978}, \cite{Arfken1971}, \cite{Gradshteyn2007}, \cite{Abramowitz1970}, \cite{Jeffrey2008}, \cite{Churchill1972}, \cite{Stalker2009}, \cite{Ash2004}) offer nil in this respect but the classic Whittaker \& Watson \cite{Whittaker1915} displays Teixeira's theorem and Laurent's series. This subject seems not to be actively studied presently. Aschenbrenner \cite{Aschenbrenner2011} has made an effort to find a power series. Other research on generating power series for analytic functions are published but are only superficially relevant regarding this work (\cite{Sobczyk2007}, \cite{Lopez2004}, \cite{Schlosser2008}). The recent research in fractional calculus Li \cite{Li2009} and Yang \cite{Yang2008} is after series expansions. Boas \cite{Boas1940} worked  with function expansions and developed some useful theorems. Campos \cite{Campos1989} has also made a good effort in this direction by using diffintegration concepts. The results are generalizations of the known series expansions. Widder \cite{Widder1928}, \cite{Widder1929} made a step towards a functional series in terms of a sequence of prescribed functions, by using Wronskian and differential operators.
\subsection{Teixeira's Theorem}
The Teixeira's theorem \cite{Teixeira1900} offers a way to expand a given function in terms of some function. It can be applied to some functions with a set of requirements. Teixeira's theorem is an extension of B\"{u}rmann's theorem. 

Using the notation of \cite{Whittaker1915}, the function $f(z)$ is analytic in region A which is ring-shaped. The outer curve is $c_1$ and the inner curve is $c_2$. $\theta(z)$ is an analytic function on and inside $c_1$ having only one simple zero at $a$ within that region and $x \in A$. For all points $z \in c_1$
\begin{equation}
\left|\theta(x)\right| < \left|\theta(z)\right|  \label{eqn5040}
\end{equation}
For all points $z \in c_2$
\begin{equation}
\left|\theta(x)\right| > \left|\theta(z)\right|  \label{eqn5050}
\end{equation}
Then we may expand the function as follows
\begin{equation}
f(x)=\sum^{\infty}_{n=0}A_n[\theta(x)]^n+\sum^{\infty}_{n=1}\frac{B_n}{[\theta(x)]^n} \label{eqn5070}
\end{equation}
where
\begin{equation}
A_n=\frac{1}{2\pi{i\cdot{n}}}\int_{c_1}\frac{f'(z)dz}{[\theta(z)]^{n}} \label{eqn5080}
\end{equation}
\begin{equation}
B_n=\frac{-1}{2\pi{i\cdot{n}}}\int_{c_2}{f'(z)[\theta(z)]^{n}}dz \label{eqn5090}
\end{equation}
This is an expansion of $f(z)$ in terms of positive and negative powers of $\theta(z)$. Due to the restrictions above, this theorem is limited in scope.
\subsection{Function Substitution}
While lacking a general method, substitutions have sometimes been used. The independent variable is replaced with a function in some known series formula, like an ordinary Taylor's power series. One is then anticipating that somehow it would sort out to some sensible expansion for his use. For example, the transcendental functions like $exp(z)$, $sin(z)$ are used and the binomial expansion too. 

Taylor's power series  for a function $f(z)$ in terms of a variable $z \in C$, may be used as a basis for substitution. 
\begin{equation}
f(z)=\sum^{\infty}_{n=0}\frac{(z-z_0)^{n}}{n!}[\frac{d^{n}f(z)}{dz^n}]_{z_0} \label{eqn10}
\end{equation}
with $f(z)$ being analytic around some point $z_0$ over the interior of a circle
\begin{equation}
r=\left|z-z_0\right| \label{eqn20}
\end{equation}
The variable $z$ is here complex in the most general analytic case. Since the bracketed term in equation (\ref{eqn10}) is a constant $b_n$, we can directly substitute $s(z)$ for $z$ to get 
\begin{equation}
f(s(z))=\sum^{\infty}_{n=0}\frac{(s(z)-s(z_0))^{n}}{n!}b_n \label{eqn30}
\end{equation}
This unfortunate method is the only practical one available in addition to the inverse-composite method. The equation (\ref{eqn30}) is not what we set out to achieve in the first place. It becomes very difficult to use if we need to expand the function $f(z)$ without the composite function $s(z)$ in the argument. 
\subsection{Problem Setup}
The principal motivation for this work is to find for a function $f(z)$ a power series in terms of a function $s(z)$, both being arbitrary. The required equation is of the following form
\begin{equation}
f(z)=\sum^{\infty}_{n=0}(s(z)-s(z_0))^{n}c_n \label{eqn40}
\end{equation}
Here $c_n$ are coefficients with no dependence on $z$.

There is a need for a simple functional power series with as few restrictions as possible on the functions. We have at least two options available to solve equation (\ref{eqn40}). We either derive an inverse-composite function pair to generate a new function $h(s)$ for which to apply (\ref{eqn30}) (see Appendix A) or we derive some new functional series. In the following section we derive the new functional power series. In the next section, some sample cases are solved exhibiting a few features of it. Appendix B introduces the remainder term for the new functional power series. The composite function's derivatives are processed in Appendix C.
\section{The Functional Power Series}
The Taylor's power series for an analytic function $k(s)$ with the variable $s=s(z) \in C$, after we mark $s_0=s(z_0)$, will be
\begin{equation}
k(s)=\sum^{\infty}_{n=0}\frac{(s-s_0)^{n}}{n!}[\frac{d^{n}k(s)}{ds^n}]_{s_0} \label{eqn80}
\end{equation}
We attempt to derive a series of the form of equation (\ref{eqn40}) from this. The equation (\ref{eqn9130}) for a composite function in Appendix C is suitable for advancing equation (\ref{eqn80}). This brings us to the simplified equation
\begin{equation}
k(s(z))=\sum^{\infty}_{n=0}\frac{(s(z)-s(z_0))^{n}}{n!}[(\frac{1}{s'(z)}\frac{d}{dz})^{n}{k(s(z))}]_{z_0} \label{eqn140}
\end{equation}
Equation (\ref{eqn140}) is not yet in the form we wanted. Since $k(s(z))$ is arbitrary we can replace it with a function $f(z)$ obtaining 
\begin{equation}
f(z)=\sum^{\infty}_{n=0}\frac{(s(z)-s(z_0))^{n}}{n!}[(\frac{1}{s'(z)}\frac{d}{dz})^{n}{f(z)}]_{z_0} \label{eqn150}
\end{equation}
The necessary but not sufficient conditions for this power series to converge, are
\begin{equation}
s'(z_0)\neq{0} \label{eqn160}
\end{equation}
and
\begin{equation}
\left|{s(z_0)}\right| < {\infty} 
\label{eqn170}
\end{equation}
Equation (\ref{eqn150}) is the functional power series for an arbitrary function $f(z)$ in terms of another arbitrary function $s(z)$. Additional parameters can be planted to $s(z)$. Also the point of focus $z_0$ is a free parameter as long as the basic conditions (\ref{eqn160}) and (\ref{eqn170}) hold and the resulting series converges. When $s(z)\rightarrow z$ this series will approach a Taylor's series. 
\section{Simple Results}
In the following we present example cases solved with the new functional power series (\ref{eqn150}), assuming $z \in C$. It seems to be typical, that the functional power series of a simple function in terms of another simple function turns out to be very simple in form. It is usually either an exponential, a logarithmic or a binomial series. 

\subsection{Rational Power Function in Terms of a Power Function and Vice Versa}
We look for the series of the function
\begin{equation}
f(z)=\frac{1}{1-2^{1-z}}         \label{eqn200}
\end{equation}
in terms of
\begin{equation}
s(z)=2^{-z}  \label{eqn220}
\end{equation}
By using the equation (\ref{eqn150}) we get
\begin{equation}
f(z)=\frac{1}{1-2^{1-z{_0}}}\sum^{\infty}_{n=0}\frac{2^n\cdot{(2^{-z}-2^{-z_0})^{n}}}{(1-2^{1-z_0})^{n}}  \label{eqn240}
\end{equation}
This equation is identified as formally being a binomial series and after sorting it out we end up with an identity. If we attempt to find the series of
\begin{equation}
f(z)=r^{-z} \label{eqn300}
\end{equation}
in terms of
\begin{equation}
s(z)=\frac{1}{1-r^{1-z}}         \label{eqn320}
\end{equation}
we obtain
\begin{equation}
f(z)=\frac{1}{r^{z_0}}-\frac{(1-r^{1-z{_0}})}{r}[\sum^{\infty}_{n=0}(1-\frac{1-r^{1-z_0}}{1-r^{1-z}})^n-1]  \label{eqn340}
\end{equation}
This equation is again a binomial series and leads to an identity. 
\subsection{Power Function in Terms of a Power Function}
The function
\begin{equation}
f(z)=k^{-z} \label{eqn400}
\end{equation}
can be expanded in terms of
\begin{equation}
s(z)=M^{-z}         \label{eqn420}
\end{equation}
resulting in 
\begin{equation}
f(z)=\frac{1}{k^{z_0}}\sum^{\infty}_{n=0}\frac{(\frac{M^{z_0}}{M^{z}}-1)^n\cdot{ln(k)\cdot{ln(\frac{k}{M})}\cdot{ln(\frac{k}{M^2})}\cdot{ln(\frac{k}{M^3})}...ln(\frac{k}{M^{n-1}})}}{n!(ln(M))^n}  \label{eqn440}
\end{equation}
This is not any common series and likely cannot be derived from a Taylor's series in a trivial way. It is notable that this series will terminate when
\begin{equation}
n=\frac{ln(k)}{ln(M)}         \label{eqn460}
\end{equation}
since $n \in N$ meaning that the radices logarithms have an integer in common. Termination is a characteristic option for the functional power series in some cases. The particular case here is a function expanded in terms of its own power $\beta$. This can be proved from the equation (\ref{eqn150}) to be a general property. Termination happens if the value of the  power $\beta$ is for some integer $n=1,2,3...$
\begin{equation}
\beta=\frac{1}{n}         \label{eqn470}
\end{equation}
The expansion of a function in terms of its own power is the following
\begin{equation}
f(z)=f(z_0)[1+\sum^{\infty}_{n=1}\frac{((\frac{f(z)}{f(z_0)})^{\beta}-1)^n}{n!\beta^n}\cdot{(1-\beta)(1-2\beta)}...(1-(n-1)\beta)]  \label{eqn480}
\end{equation}

\subsection{Rational Function in Terms of a Rational Function}
We need the series of the rational function
\begin{equation}
f(z)=\frac{1}{(z-a)(z-a)}         \label{eqn1200}
\end{equation}
in terms of
\begin{equation}
s(z)=\frac{1}{z-a}  \label{eqn1220}
\end{equation}
By applying the equation (\ref{eqn150}) again we obtain
\begin{equation}
f(z)=\frac{-1}{a^2}-\frac{2}{a(z-a)}+(\frac{1}{z-a}+\frac{1}{a})^2  \label{eqn1240}
\end{equation}
This is a partial fraction "`solution"' to the famous degenerate case. This is an example of a series of a function in terms of its power leading to a terminating series, becoming an identity.
\subsection{Rational Function in Terms of a Sinus Function, an Approximation}
We want to approximate the function $f(z)$ around origin
\begin{equation}
f(z)=\frac{1}{1+z}         \label{eqn1300}
\end{equation}
with a few terms of
\begin{equation}
s(z)=sin(z)  \label{eqn1320}
\end{equation}
The equation (\ref{eqn150}) gives the first four terms
\begin{equation}
f(z)\approx{1-sin(z)+sin^2(z)-\frac{7}{6}sin^3(z)}  \label{eqn1340}
\end{equation}

\begin{figure}
	\centering
		\includegraphics[width=0.96\textwidth]{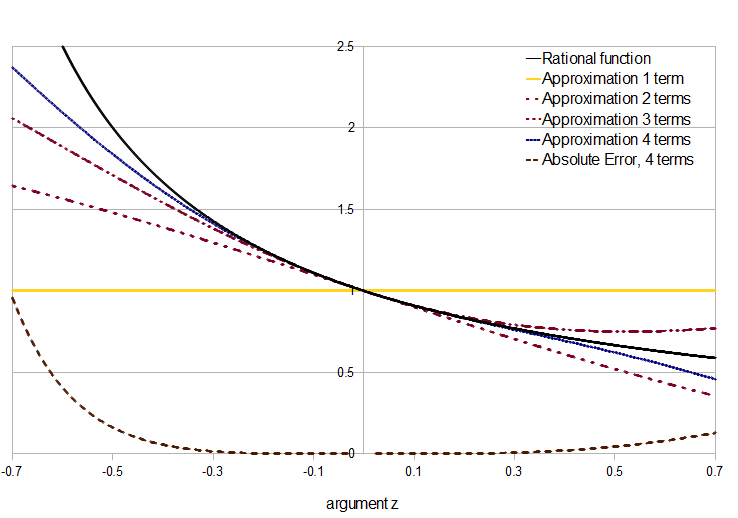}
	\caption{Approximation of the elementary rational function in terms of the sinus function. A few first approximations are shown and with real argument values only}
	\label{fig:sinx}
\end{figure}

In Fig. \ref{fig:sinx} we have a real variable case. We can observe that even by using just a few terms will give a good approximation, although the sinus function was not particularly selected to match the rational function at the origin.
\section{Discussion}
The equation (\ref{eqn150}) represents a new functional power series for an arbitrary function in terms of another arbitrary function. Formally, there are no limits in this respect, except the conditions (\ref{eqn160}) and (\ref{eqn170}) must be valid at the point of focus. The functional series can be applied to expanding any function in terms of any function as long as the derivatives
\begin{equation}
\frac{d^{n}s(z)}{dz^n}    \label{eqn7160}
\end{equation}
and 
\begin{equation}
\frac{d^{n}f(z)}{dz^n} \label{eqn7180}
\end{equation}
exist for all $n$. Analyticity, of course, will guarantee this. 

Partial (truncated) series can be used in complicated cases where the pattern of the coefficients $c_n$ is unclear. The resulting error can be estimated with a remainder term. The functional series requires only elementary differentiation and evaluation at $z_0$. Thus it is suitable for hand calculations and symbolic machine processing. The $s(z)$ function may contain free parameters. 

Termination of the series is characteristic in some cases. Terminating series create identities which may be very complicated. 

The functional series is useful for approximating the behavior of a function around a point of interest in terms of some other function. The function to be used can be selected to possibly speed up convergence and minimize the number of terms for the required accuracy. Asymptotic behavior of a function may be studied with the new series. 

As a byproduct of this work, the equations (\ref{eqn9130}) and (\ref{eqn9132}) in Appendix C are new expressions for the n'th derivative of a composite function. 


\appendix
\section{Appendix. The Inverse-Composite Function and The Taylor's Series of Functions}
As noted above, the Taylor's series (refer to equation (\ref{eqn30})) can be used as basis for a functional series. Intuitively, we can create an inverse-composite function pair as follows. The composite function is
\begin{equation}
s=s(z)
\end{equation}
and the inverse would be, if solvable
\begin{equation}
z=g(s)=s^{-1}(s(z)) \label{eqn7550}
\end{equation}
Thus
\begin{equation}
f(z)=f(g(s(z)))=h(s(z)) \label{eqn7560}
\end{equation}
We can use $h(s(z))$ instead of $f(z)$ and use the Taylor's series of $h(s)$, assuming it exists, to yield
\begin{equation}
f(z)=h(s(z))=\sum^{\infty}_{n=0}\frac{(s(z)-s(z_0))^{n}}{n!}[\frac{d^{n}h(s)}{ds^n}]_{s_0} \label{eqn7570}
\end{equation}
To have any practical use of this equation one needs to be able to solve in algebraic form the inverse function $g(s)$. This is possible only with a small set of functions. 

\section{Appendix. Remainder Term for the Functional Power Series}
The remainder term of a series expansion equals the error caused by breaking the series after the N'th term. 
Assuming $f(z)$ and $s(z)$ are analytic (referring to equation (\ref{eqn150})) and by using equation (\ref{eqn9130}) in Appendix C, it is not difficult to prove by induction the expression
\begin{equation}
f(z)=\sum^{N}_{n=0}\frac{(s(z)-s(z_0))^{n}}{n!}[(\frac{1}{s'(z)}\frac{d}{dz})^{n}{f(z)}]_{z_0}+   \nonumber
\end{equation}
\begin{equation}
\frac{1}{N!}\int^{s(z)}_{s(z_0)}d\chi \cdot(s-\chi)^{N}\left[(\frac{1}{\frac{d\chi}{dz}}\frac{d}{dz})^{N+1}f(z)\right]_{z=z^{-1}(\chi)} \label{eqn8010}
\end{equation}
$z=z^{-1}(\chi)$ is the inverse function of $s(z)$. For the real variable case, by using the mean value theorem and integrating, we get the latter part which is the remainder
\begin{equation}
R_N(s)=\frac{(s(z)-s(z_0))^{N+1}}{(N+1)!}\left[(\frac{1}{\frac{ds(z)}{dz}}\frac{d}{dz})^{N+1}f(z)\right]_{z=z^{-1}(\zeta)} \label{eqn8020}
\end{equation}
$\zeta$ is some value of $s$ between $s(z_0)$ and $s(z)$. This result is analogous to the Lagrange remainder for real functions.
For the complex variable case, by using the corresponding mean value theorem IV by Curtiss \cite{Curtiss1907} we get after integration
\begin{equation}
R_N(s)=\frac{[(\theta-1)(s(z)-s(z_0))]^{N+1}}{(N+1)!}\left[(\frac{1}{\frac{ds(z)}{dz}}\frac{d}{dz})^{N+1}f(z)\right]_{z=z^{-1}(s(z_0))} \label{eqn8040}
\end{equation}
being actually equal to
\begin{equation}
R_N(s)=\frac{[(\theta-1)(s(z)-s(z_0))]^{N+1}}{(N+1)!}\left[(\frac{1}{\frac{ds(z)}{dz}}\frac{d}{dz})^{N+1}f(z)\right]_{z_0} \label{eqn8050}
\end{equation}
Here $\theta$ has some value according to
\begin{equation}
\left|{\theta-1}\right|<1   \label{eqn8070}
\end{equation}

\section{Appendix. Derivatives of a Composite Function}

The derivatives of a composite function has been studied earlier and usually combinatoric methods are used (refer to \cite{Riordan1946}, \cite{McKiernan1956}, \cite{Floater2007}). The resulting expressions are complicated and not suitable for our equation (\ref{eqn80}).

We need to solve the n'th derivative of a composite function
\begin{equation}
\frac{d^{n}k(s)}{ds^n} \label{eqn9000}
\end{equation}
with the function $s(z)$ as the composite. An expression with derivatives by $z$ is required. We take the variable $z \in C$ as an independent variable and unravel all the derivative terms. The first derivative would be
\begin{equation}
\frac{dk}{ds}=\frac{1}{s'(z)}\frac{dk}{dz} \label{eqn9090}
\end{equation}
The second derivative is
\begin{equation}
\frac{d^2{k}}{ds^2}=\frac{d}{dz}(\frac{1}{s'(z)}\frac{dk}{dz})\cdot\frac{dz}{ds} \label{eqn9100}
\end{equation}
If differentiation is performed fully, we end up with a slightly complicated expression. It is obvious that continuing in this way will lead to an enormous complexity. This is very difficult to be used any further as such. In order to circumvent the generation of progressively complicated terms in higher derivatives, we proceed as follows. Equation (\ref{eqn9100}) becomes, without fully completing the differentiation,
\begin{equation}
\frac{d^2{k}}{ds^2}=\frac{1}{s'(z)}\cdot\frac{d}{dz}(\frac{1}{s'(z)}\frac{dk}{dz}) \label{eqn9110}
\end{equation}
Continuing along this route leads to the n'th derivative
\begin{equation}
\frac{d^n{k}}{ds^n}=\frac{1}{s'(z)}(\frac{d}{dz}(\frac{1}{s'(z)}(\frac{d}{dz}(\frac{1}{s'(z)}(\frac{d}{dz}\cdot\cdot\cdot\cdot(\frac{dk}{dz})))))) \label{eqn9120}
\end{equation}
We can rearrange the brackets to obtain 
\begin{equation}
\frac{d^n{k(s)}}{ds^{n}}=(\frac{1}{s'(z)}\frac{d}{dz})^{n}{k(s(z))} \label{eqn9130}
\end{equation}
This is the required end result.

We can also perform analogous steps from equation (\ref{eqn9090}) in the opposite direction to get
\begin{equation}
\frac{d^n{k(s(z))}}{dz^{n}}=(s'(z)\frac{d}{ds})^{n}{k(s)} \label{eqn9132}
\end{equation}

\end{document}